\documentstyle[12pt]{article}
\newcommand{\height}{{\mathrm{height}\,}}
\newcommand{\Spec}{{\mathrm{Spec}\,}}
\newcommand{\Frac}{{\mathrm{Frac}\,}}
\newtheorem{theorem}{Theorem}[section]
\newtheorem{lemma}[theorem]{Lemma}

\newtheorem{corollary}[theorem]{Corollary}
\newtheorem{remar}[theorem]{Remark}
\newenvironment{proof}{Proof:\ \ \ }{\QED}
\newenvironment{remark}{\begin{remar}\rm}{\end{remar}}

\newcommand{\QED}{{\unskip\nobreak\hfil\penalty50%
\hskip1em\hbox{}\nobreak\hfil $\Box$%
\parfillskip=0pt \finalhyphendemerits=0 \par\medskip\noindent}}
\newcommand{\bfind}[1]{\index{#1}{\bf #1}}

\newcommand{\n}{\par\noindent}

\newcommand{\sn}{\par\smallskip\noindent}
\newcommand{\mn}{\par\medskip\noindent}

\newcommand{\pars}{\par\smallskip}
\newcommand{\parm}{\par\medskip}

\newcommand{\sep}{^{\rm sep}}

\newcommand{\chara}{\mbox{\rm char}\,}
\newcommand{\trdeg}{\mbox{\rm trdeg}\,}

\newcommand{\adresse}{\par\bigskip \small\rm
 Fraunhofer Institut Techno- und Wirtschaftsmathematik,\par
 Gottlieb-Daimler-Strasse, Geb.\ 49, D--67663 Kaiserslautern,
 Germany\par
 email: knaf@itwm.fhg.de
  \par\bigskip
 Mathematical Sciences Group, %\par
 University of Saskatchewan, \par
 106 Wiggins Road, %\par
 Saskatoon, Saskatchewan, Canada S7N 5E6 \par
 email: fvk@math.usask.ca \ \ --- \ \ home page:
 http://math.usask.ca/$\,\tilde{ }\,$fvk/index.html
}
\font\tenlv=msbm10  scaled 1200
\font\sevenlv=msbm7 scaled 1200
\font\fivelv=msbm5  scaled 1200

\def\lv #1{{\mathchoice{{\hbox{\tenlv #1}}}{{\hbox{\tenlv #1}}}
{{\hbox{\sevenlv #1}}}{{\hbox{\fivelv #1}}}}}
\newcommand{\N}{\lv N}
\newcommand{\Q}{\lv Q}
\newcommand{\R}{\lv R}
\newcommand{\Z}{\lv Z}

\newcommand{\A}{\lv A}

\begin{document}
\title{Abhyankar places admit local uniformization in any
characteristic}
\author{Hagen Knaf and Franz--Viktor Kuhlmann\footnote{
The second author would like to thank Hans Schoutens, Peter Roquette,
Dale Cutkosky and Olivier Piltant for many inspiring discussions.}}
\date{19.\ 3.\ 2003}
\maketitle
\begin{abstract}\noindent
We prove that every place $P$ of an algebraic function field $F|K$ of
arbitrary characteristic admits local uniformization, provided that the
sum of the rational rank of its value group and the transcendence degree
of its residue field $FP$ over $K$ is equal to the transcendence degree
of $F|K$, and the extension $FP|K$ is separable. We generalize this
result to the case where $P$ dominates a regular local Nagata ring
$R\subseteq K$ of Krull dimension $\dim R\leq 2$, assuming that the
valued field $(K,v_{P})$ is defectless, the factor group $v_{P}F/v_{P}K$
is torsion-free and the extension of residue fields $FP|KP$ is
separable. The results also include a form of monomialization.
Further, we show that in both cases, finitely many Abhyankar places
admit simultaneous local uniformization on an affine scheme if they
have value groups isomorphic over $v_{P}K$.
\end{abstract}
%
%---------------------------------------------------------------
%
\section{Introduction and main results}
In [Z], Zariski proved the Local Uniformization Theorem for places of
algebraic function fields over base fields of characteristic 0. In [Z3],
he uses this theorem to prove resolution of singularities for algebraic
surfaces in characteristic 0. This result was generalized by Abhyankar
to the case of positive characteristic [A4] and to the case of
arithmetic surfaces over a Dedekind domain [A5]. More recently de Jong
[dJ] proved that the singularities of an algebraic or arithmetic variety
$X$ can be resolved by successively applying a finite number of
morphisms called alterations. An alteration $f:Y\rightarrow X$ is a
composition $f=g\circ h$ of a birational isomorphism
$h:X^\prime\rightarrow X$ and a finite morphism $g:Y\rightarrow
X^\prime$. In general it leads to a finite extension $K(Y)|K(X)$ of
function fields. Applying de Jong's results to a proper algebraic
variety over a field of positive characteristic or to a proper
arithmetic variety over a discrete valuation ring implies Local
Uniformization in these cases, provided one allows finite extensions of
the variety's function field or of the discrete valuation ring,
respectively. As the resolution of singularities in the birational sense
is still an open problem, one is interested in generalizations of the
Local Uniformization Theorem.

 %Inserted by TeXtelmExtel

In the present paper we investigate the class of so-called Abhyankar
places of a function field $F|K$. We prove that Abhyankar places that
are trivial on $K$ admit local uniformization on algebraic varieties
without any extension of $F$, provided they have a
separable residue field extension. We also consider the case of
Abhyankar places dominating a regular local ring $R\subset K$ and
provide sufficient conditions for local uniformization on $R$-models of
$F|K$. The core of our method is purely of valuation-theoretic nature.
Before stating the main results precisely, we have to introduce
the necessary terminology.

 %Inserted by TeXtelmExtel

\smallskip
Throughout this paper the term \textit{function field} will always
mean \textit{algebraic function field}. We consider places $P$ of a
function field $F|K$ not necessarily inducing the identity on the
constant field $K$. Places that do induce the identity on $K$ (or an
isomorphism) are called \textbf{$K$-trivial}. The valuation associated
with $P$ will be denoted by $v_P^{ }\,$, the $v_P^{ }$-value of an
element $a$ by $v_P^{ }a$ and consequently, the value group of $v_P$ on
$F$ by $v_P^{ }F$. Places are considered to operate on the right: The
residue of an element $a\in F$ is denoted by $aP$ and consequently, $FP$
stands for the residue field of $F$ with respect to $P$. We
frequently do not distinguish between a place $P$ on the field $F$ and
its restrictions to subfields $E\subseteq F$. Since we are usually
working with one fixed place this does not lead to confusion. The
valuation ring of $P$ on a subfield $E\subseteq F$ is denoted by
${\cal O}_E$. If one has to distinguish between several places on the
same field $F$ the notation ${\cal O}_P$ is used for the valuation ring
of $P$ on $F$. The maximal ideal of the valuation rings ${\cal O}_{P}$
and ${\cal O}_{E}$ is denoted by ${\cal M}_{P}$ and ${\cal M}_{E}$,
respectively. Finally it should be mentioned that by abuse of language
we refer to a pair $(F,P)$ consisting of a field $F$ and a place $P$ of
$F$ as a valued field, and to a pair $(F|K,P)$, $F|K$ an extension of
fields, as an extension of valued fields, or as a valued function field
if $F|K$ is a function field.

 %Inserted by TeXtelmExtel

For every place $P$ of a function field $F|K$, we have the following
inequality:
\begin{equation}                            \label{Abhie}
\trdeg F|K \>\geq\> \trdeg FP|KP \,+\, \dim_\Q ((v_P^{ }F/v_P{ }K)
\otimes\Q )\;.
\end{equation}
Note that $\dim_\Q (v_P^{ }F/v_P{ }K\otimes\Q )$ is the \bfind{rational
rank} of the abelian group $v_P^{ }F/v_P{ }K$, i.e., the maximal number
of rationally independent elements in $v_P^{ }F/v_P{ }K$. We call $P$ an
\bfind{Abhyankar place of $F|K$} if equality holds in (\ref{Abhie}).

 %Inserted by TeXtelmExtel

In the context of local uniformization we shall be concerned with the
prime factorization in a regular local ring and use the following terminology:
Let ${\cal O}$ be a commutative ring and $H\subseteq {\cal O}$. An
element $a\in{\cal O}$ is called an \textbf{${\cal O}$-monomial in $H$} if
\[
a=u\prod\limits_{i=1}^d h_i^{\mu_i},\; u\in{\cal O}^\times ,\; h_i\in H,\;
\mu_i\in\N_0 ,\;  i=1,\ldots ,d,
\]
holds, where $\N_0:=\N\cup \{0\}$.

 %Inserted by TeXtelmExtel

We investigate local uniformization over regular base rings: Let
$(F|K,P)$ be a valued function field, and let $R\subseteq {\cal O}_{P}$
be a subring of $K$ with $\Frac R=K$. Moreover, let $Z\subset
{\cal O}_{P}$ be a finite set. The pair $(P,Z)$ is called
\textbf{$R$-uniformizable} if there exists an integral separated
$R$-scheme $X$ of finite type with field of rational functions $F$ ---
an \textbf{$R$-model of $F$} --- such that $P$ is centered in a regular
point $x\in X$ and $Z$ is a subset of the local ring ${\cal O}_{X,x}$ at
$x$. We also use the phrase \textbf{$R$-uniformizable on $X$} if this
situation is present. The place $P$ is called $R$-uniformizable iff
$(P,\emptyset )$ is $R$-uniformizable. Note that including the finite
set $Z$ in the problem of uniformization allows to prove statements
like: Given an integral $R$-scheme $X$ such that $P$ is centered in the
(singular) point $x\in X$ there exists an integral $R$-scheme $Y$ and a
morphism $\pi : U\rightarrow V$, where $U\subseteq Y$ and $V\subseteq X$
are open subschemes such that $x\in V$ and $P$ is centered in a regular
point $y\in\pi ^{-1}x$.

 %Inserted by TeXtelmExtel

Throughout the paper we restrict ourselves to the case of a {\em regular
local ring} $R$ with maximal ideal $M={\cal M}_{P}\cap R$, i.e., we
assume that ${\cal O}_{P}$ dominates $R$.

 %Inserted by TeXtelmExtel

We prove a Local Uniformization Theorem for $K$-trivial Abhyankar
places in arbitrary characteristic:
\begin{theorem}                                \label{MT}
Let $P$ be a $K$-trivial Abhyankar place of the function field $F|K$,
and assume that $FP|K$ is separable. Take a finite set $Z\subset
{\cal O}_{P}\,$. Then the pair $(P,Z)$ is $K$-uniformizable on a variety
$X$ such that $P$ is centered in a smooth point $x\in X$ and $\dim
{\cal O}_{X,x}=\dim_{\Q} (v_{P}F\otimes\Q)$. Moreover, $X$ can be
chosen such that all $\zeta\in Z$ are ${\cal O}_{X,x}$-monomials in
$\{a_1,\ldots ,a_d\}$ for some regular parameter system $(a_1,\ldots ,a_d )$
of ${\cal O}_{X,x}$.
\end{theorem}

 %Inserted by TeXtelmExtel

Theorem~\ref{MT} essentially is proved by embedding $F$ in the
field of fractions of the strict henselization of the valuation ring
${\cal O}_{K(T)}= {\cal O}_P\cap K(T)$ for a suitable transcendence
basis $T\subset F$ of $F|K$. The methods used in the proof of
Theorem~\ref{MT} are applicable even if $P$ is not trivial on $K$. They
then lead to a uniformization result for integral schemes of finite type
over certain base rings $R\subset K$. For the definition of the notion
``defectless'', see Section~\ref{sect3}.

 %Inserted by TeXtelmExtel

\begin{theorem}                                \label{MT2}
Let $P$ be an Abhyankar place of the function field $F|K$, which is
non-trival on $K$. Assume that $(K,P)$ is defectless, $FP|KP$ is
separable and the group $v_{P}F/v_{P}K$ is torsion-free.

 %Inserted by TeXtelmExtel

Let $R\subset K\cap {\cal O}_{P}$, $\Frac R=K$, be a noetherian, regular local
ring with maximal ideal $M={\cal M}_{P}\cap R$ and of dimension $\dim R\leq 2$.
Assume that $R$ is a Nagata ring if $\dim R=2$.

 %Inserted by TeXtelmExtel

Then for every finite set $Z\subset {\cal O}_{P}$ the pair $(P,Z)$ is
$R$-uniformizable on an $R$-scheme $X$ such that the center $x\in X$ of
$P$ on $X$ satisfies:
\begin{itemize}
\item  $\dim {\cal O}_{X,x}=\dim_{\Q}(v_{P}F/v_{P}K\otimes\Q)+1$ if
${\cal O}_{K}$ is a discrete valuation ring.
\item  $\dim {\cal O}_{X,x}=\dim_{\Q}(v_{P}F/v_{P}K\otimes\Q)+2$ in the
remaining cases.
\end{itemize}
Moreover, $X$ can be
chosen such that all $\zeta\in Z$ are ${\cal O}_{X,x}$-monomials in
$\{a_1,\ldots ,a_d\}$ for some regular parameter system $(a_1,\ldots ,a_d )$
of ${\cal O}_{X,x}$.
\end{theorem}

 %Inserted by TeXtelmExtel

Some remarks concerning the condition $\dim R\leq 2$ may be helpful at
this point. Recall that the domain $R$ is called \textbf{Nagata}, if the
integral closure of every factor ring $R/p$, $p\in\Spec R$, in every
finite extension of $\Frac R/p$ is finite. In the case of $\dim R=1$,
the ring $R$ is a discrete valuation ring and moreover, in the situation
of the theorem, we have $R={\cal O}_K$. It is well-known that then the
property of being Nagata is equivalent to the valued field $(K,P)$ being
defectless. {\em We conclude that all base rings appearing in
Theorem~\ref{MT2} are Nagata.}

 %Inserted by TeXtelmExtel

A further important ring-theoretic notion that we will have to use is
universal catenarity: The domain $R$ is called \textbf{universally
catenary} if every polynomial ring $R[X_1,\ldots ,X_n]$, $n\in\N$, has
the property: for every pair of prime ideals $p,q\in R[X_1,\ldots ,X_n]$
with $p\subset q$, all non-refinable chains of primes $p=:p_0\subset
p_1\subset\ldots \subset p_\ell :=q$ have a common finite length $\ell$
(depending on $p,q$).

 %Inserted by TeXtelmExtel

Examples of noetherian, universally catenary rings are the
Cohen-Macaulay rings; in particular any 2-dimensional, normal,
noetherian ring is universally catenary. Furthermore, it is known that
valuation domains are universally catenary [Nag], whether they are
noetherian or not. {\em We conclude that all base rings appearing
in Theorem~\ref{MT2} are universally catenary.}

 %Inserted by TeXtelmExtel

Every universally catenary domain $R$ satisfies the \textbf{altitude
formula}: for every prime $q\in\Spec A$ of every domain $A$ finitely
generated over $R$ the equation
\begin{equation}                              \label{altitude}
\height q +\trdeg (A/q |R/p)=\height p+\trdeg (A|R),\; p:=q\cap R,
\end{equation}
holds.

 %Inserted by TeXtelmExtel

The condition $\dim R\leq 2$ can be replaced by a more general but
rather technical condition that ensures the existence of certain
monoidal transforms of $R$ along the valuation $v_{P}$.

 %Inserted by TeXtelmExtel

\parm
In [K4] it is shown that the $K$-trivial Abhyankar places lie dense in
the Zariski space of all $K$-trivial places of $F|K$, with respect to a
``Zariski patch topology''. This topology is finer than the Zariski
topology (but still compact); its basic open sets are the sets of the
form
\[\{P\mid \mbox{$P$ a place of $F|K$ such that }a_1 P\ne 0,\ldots,
a_k P\ne 0\,;\,b_1 P=0,\ldots,b_\ell P=0\}\]
with $a_1 ,\ldots,a_k,b_1,\ldots, b_\ell\in F\setminus \{0\}$.
Theorem~\ref{MT} thus yields:

 %Inserted by TeXtelmExtel

\begin{corollary}
The $K$-uniformizable places of $F|K$ lie dense in the Zariski space
of~$F|K$, with respect to the Zariski patch topology, provided that
$K$ is perfect.
\end{corollary}
We conclude with a generalization of Theorem~\ref{MT2}. It gives
simultaneous uniformization on an affine scheme for finitely many
Abhyankar places which satisfy the assumptions of Theorem~\ref{MT2} and
have isomorphic value groups.

 %Inserted by TeXtelmExtel

\begin{theorem}                                       \label{MT3}
Let $F|K$ be a function field equipped with Abhyankar places $P_1,\ldots
,P_\ell$ having a common not necessarily trival restriction $P$ to $K$.
Assume that $(K,P)$ is defectless, that for $i=1,\ldots ,\ell$ the
extensions $FP_i|KP$ are separable, that there exist isomorphisms
$\phi_i: v_{P_i}F\rightarrow v_{P_1}F$ of totally ordered abelian groups
such that $\phi_i |_{v_{P}K}=\textrm{id}$, and that the
group $v_{P_1}F/v_{P}K$ is torsion-free.

 %Inserted by TeXtelmExtel

Further let $R\subseteq {\cal O}_{P}$, $\Frac R=K$, be a noetherian, regular local ring with
$\dim R\leq 2$ and maximal ideal $M={\cal M}_{P}\cap R$. Assume that $R$
is a Nagata ring if $\dim R=2$.

 %Inserted by TeXtelmExtel

Then for every finite set $Z\subset\cap_{i=1}^\ell {\cal O}_{P_i}$ there
exists a morphism $\pi : X\rightarrow Y$ of affine $R$-schemes, where $X$
is an $R$-model of $F$, $Y$ has rational function field $K(Y)$
and $F|K(Y)$ is finite, such that the places $P_1,\ldots ,P_\ell$
are centered in regular points $x_1,\ldots ,x_\ell\in X$ with
$Z\subset \cap_{i=1}^\ell {\cal O}_{X,x_i}$, and $\pi x_1,\ldots
,\pi x_\ell\in Y$ are regular on $Y$.

 %Inserted by TeXtelmExtel

If $R=K$, then the morphism $\pi$ can be chosen such that the points
$x_1,\ldots ,x_\ell$ and their images are smooth over $K$.
\end{theorem}
An interesting application of Theorem~\ref{MT3} is based on the following
valuation-theoretic fact proved in [K3]: Let $A$ be a domain finitely generated
over the field $K$ and let $q_1,\ldots ,q_\ell\in\Spec A$ be primes with
$d=\dim A/q_i$, $i=1,\ldots ,\ell$. Then there exist Abhyankar places $P_1,\ldots ,P_\ell$
of $F=\Frac A|K$ possessing isomorphic value groups. A straightforward consequence
now is:
\begin{corollary}
Let $X=\Spec A$ be an affine variety over the perfect field $K$. Let
$x_1,\ldots ,x_\ell$ $\in X$ be points of equal codimension. Then there
exist Abhyankar places $P_1,\ldots ,P_\ell$ of the function field
$K(X)|K$ and an affine model $Y$ of $K(X)$ such that:
\begin{itemize}
\item $P_i$ is centered in $x_i$ for $i=1,\ldots ,\ell$.
\item The centers $y_1,\ldots ,y_\ell$ of $P_1,\ldots ,P_\ell$ on $Y$ are smooth.
\item $A\subset {\cal O}_{Y,y_i}$ for $i=1,\ldots ,\ell$.
\end{itemize}
\end{corollary}

 %Inserted by TeXtelmExtel

%
%---------------------------------------------------------------
%
\section{Valuation independence}
The following theorem, together with Theorem~\ref{ai} below, gives the
motivation for the definition of the distinguished class of Abhyankar
places. For its proof see [B], Chapter VI, \S10.3, Theorem~1.

 %Inserted by TeXtelmExtel

\begin{theorem}                                \label{prelBour}
Let $(F|K,P)$ be an extension of valued fields. Take elements $x_i,y_j
\in F$, $i\in I$, $j\in J$, such that the values $v_P^{ }x_i\,$, $i\in
I$, are rationally independent over $v_P^{ }K$, and the residues $y_jP$,
$j\in J$, are algebraically independent over $KP$. Then the elements
$x_i,y_j$, $i\in I$, $j\in J$, are algebraically independent over $K$.

 %Inserted by TeXtelmExtel

Moreover, if we write
\[f\>=\> \displaystyle\sum_{k}^{} c_{k}\,
\prod_{i\in I}^{} x_i^{\mu_{k,i}} \prod_{j\in J}^{} y_j^{\nu_{k,j}}\in
K[x_i,y_j\mid i\in I,j\in J]\]
in such a way that for every $k\ne\ell$
there is some $i$ s.t.\ $\mu_{k,i}\ne\mu_{\ell,i}$ or some $j$ s.t.\
$\nu_{k,j}\ne\nu_{\ell,j}\,$, then
\[v_P^{ }f\>=\>\min_k\left(v_P^{ }\, (c_k \prod_{i\in I}^{}
x_i^{\mu_{k,i}}\prod_{j\in J}^{} y_j^{\nu_{k,j}})\right)\>=\>
\min_k\left(v_P^{ }\,c_k\,+\,\sum_{i\in I}^{} \mu_{k,i} v_P^{ }
x_i\right)\;.\]
That is, the value of the polynomial $f$ is equal to the least of the
values of its monomials. In particular, this implies:
\begin{eqnarray*}
v_P^{ }K(x_i,y_j\mid i\in I,j\in J) & = & v_P^{ }K\oplus
\bigoplus_{i\in I} \Z v_P^{ }x_i\\
K(x_i,y_j\mid i\in I,j\in J)P & = & KP\,(y_jP\mid j\in J)\;.
\end{eqnarray*}
It also implies that the valuation $v_P^{ }$ and the place $P$ on
$K(x_i,y_j\mid i\in I,j\in J)$ are uniquely determined by their
restrictions to $K$, the values $v_P^{ }x_i$ and the residues $y_jP$.
\end{theorem}

 %Inserted by TeXtelmExtel

Every finite extension $L$ of the valued field $(K,P)$
satisfies the {\bf fundamental inequality} (cf.\ [En]):
\begin{equation}                             \label{fundineq}
[L:K]\>\geq\>\sum_{i=1}^{\rm g} {\rm e}_i {\rm f}_i
\end{equation}
where $P_1,\ldots,P_{\rm g}$ are the distinct extensions of $P$ from $K$
to $L$, ${\rm e}_i=(v_{P_i}L:v_P^{ }K)$ are the respective ramification
indices and ${\rm f}_i=[LP_i:KP]$ are the respective inertia degrees.
Note that ${\rm g}=1$ if $(K,P)$ is henselian.
\begin{corollary}                              \label{fingentb}
Let $(F|K,P)$ be an extension of valued fields of finite transcendence
degree. Then the following inequality holds:
\begin{equation}                            \label{wtdgeq}
\trdeg F|K \>\geq\> \trdeg FP|KP \,+\, \dim_\Q ((v_P^{ }F/v_P^{ }K)
\otimes\Q )\;.
\end{equation}
If in addition $F|K$ is a function field, and if equality holds in
(\ref{wtdgeq}), then the extensions $v_P^{ }F| v_P^{ }K$ and $FP|KP$
are finitely generated. In particular, if $P$ is trivial on $K$, then
$v_P^{ }F$ is a product of finitely many copies of $\Z$, and $FP$ is
again a function field over $K$.
\end{corollary}
\begin{proof}
Choose elements $x_1,\ldots,x_{\rho},y_1,\ldots,y_{\tau}\in F$ such
that the values $v_P^{ }x_1,\ldots,v_P^{ }x_{\rho}$ are rationally
independent over $v_P^{ }K$ and the residues $y_1P,\ldots,y_{\tau} P$
are algebraically independent over $KP$. Then by the foregoing theorem,
$\rho+\tau\leq\trdeg F|K$. This proves that $\trdeg FP|KP$ and the
rational rank of $v_P^{ }F/v_P^{ }K$ are finite. Therefore, we may
choose the elements $x_i,y_j$ such that $\tau=\trdeg FP|KP$ and
$\rho=\dim_\Q ((v_P^{ }F/v_P^{ }K)\otimes\Q )$ to obtain inequality
(\ref{wtdgeq}).

 %Inserted by TeXtelmExtel

Assume that this is an equality. This means that for $F_0:=K(x_1,\ldots,
x_{\rho},y_1,\ldots,y_{\tau})$, the extension $F|F_0$ is algebraic.
Since $F|K$ is finitely generated, it follows that $F|F_0$ is finite. By
the fundamental inequality (\ref{fundineq}), this yields that $v_P^{ }
F|v_P^{ }F_0$ and $FP| F_0P$ are finite extensions. Since already
$v_P^{ } F_0|v_P^{ } K$ and $F_0P|KP$ are finitely generated by the
foregoing theorem, it follows that also $v_P^{ }F|v_P^{ }K$ and $FP|KP$
are finitely generated.
\end{proof}
\n
If equality holds in (\ref{wtdgeq}) we will either say that
\textbf{$(F|K,P)$ is without transcendence defect} or as already defined
earlier that $P$ is an Abhyankar place of $F|K$.

 %Inserted by TeXtelmExtel

%
%---------------------------------------------------------------
%
\section{Inertially generated function fields}         \label{sect3}
In this section we provide the valuation-theoretic core of the present
paper. It is a generalization of the ``Grauert--Remmert Stability
Theorem'' which is proved in [K1] (cf.\ also [K5]). To state it, we
introduce a fundamental notion: a valued field $(K,P)$ is called
\bfind{defectless} (or \bfind{stable}) if equality holds in the
fundamental inequality (\ref{fundineq}) for every finite extension
$L|K$. If $\chara KP=0$, then $(K,P)$ is defectless (this is a
consequence of the ``Lemma of Ostrowski'', cf.\ [En], [R]).

 %Inserted by TeXtelmExtel

%The \bfind{Stability Theorem} deals with the behavior of defectlessness
%in a valued function field $(F|K,P)$.
\begin{theorem}                           \label{ai}
{\bf (Generalized Stability Theorem)} \
Let $(F|K,P)$ be a valued function field without transcendence defect.
If $(K,P)$ is a defectless field, then also $(F,P)$ is a defectless
field.
\end{theorem}

 %Inserted by TeXtelmExtel

In what follows we consider concepts like the henselization of a valued
field $(K,P)$, where one has to fix an extension of $P$ and $v_{P}$ to
the algebraic closure of $K$. Thus, whenever we talk of a valued field
$(K,P)$, we will implicitly assume the valuation $v_{P}$ and the place
$P$ to be extended to the algebraic closure of $K$, the extensions
denoted by $v_{P}$ and $P$ again. Therefore, we will talk of {\it the}
henselization $K^h$, and of {\it the} {\bf absolute inertia field $K^i$
of $K$}, which we define to be the inertia field of the normal separable
extension $K\sep|K$ with respect to the given valuation $v_{P}$; here,
$K\sep$ denotes the separable-algebraic closure of $K$.

 %Inserted by TeXtelmExtel

The following lemma is proved in [K7] (and partially also in [En]):
\begin{lemma}                               \label{hdl}
A valued field $(K,P)$ is defectless if and only if its henselization
$(K^h,P)$ is.
\end{lemma}

 %Inserted by TeXtelmExtel

An extension $(L|K,P)$ of valued fields is called \bfind{immediate}
if the canonical embeddings $KP\rightarrow LP$ and $v_P^{ }K\rightarrow
v_P^{ }L$ are onto.

 %Inserted by TeXtelmExtel

\begin{corollary}                           \label{corhdie}
If $(K,P)$ is defectless, then $(K^h,P)$ does not admit proper immediate
algebraic extensions.
\end{corollary}
\begin{proof}
If $(K,P)$ is defectless, then so is $(K^h,P)$, by the foregoing lemma.
Suppose that $(L|K^h,P)$ is a finite immediate algebraic extension.
Hence, $(v_P^{ }L:v_P^{ }K^h)=1=[LP:K^hP]$. Since $(K^h,P)$ is a
henselian field, there is a unique extension of $v_P^{ }$ from $K^h$ to
$L$. Since $(K^h,P)$ is defectless, we have that $[L:K^h]=(v_P^{ }L:
v_P^{ }K^h)[LP:K^hP]=1$, showing that $L=K^h$. As every proper immediate
algebraic extension would contain a proper immediate finite extension,
it follows that $(K^h,P)$ does not admit any proper immediate algebraic
extension.
\end{proof}

 %Inserted by TeXtelmExtel

From these facts, we deduce:
\begin{theorem}                             \label{hrwtd}
Assume that $(F|K,P)$ is a valued function field without transcendence
defect such that $FP|KP$ is a separable extension, $(K,P)$ is a
defectless field and $v_{P}F/v_{P}K$ is torsion-free. Then $(F|K,P)$ is
\bfind{inertially generated}, by which we mean that there is a
transcendence basis
$T=\{x_1,\ldots,x_\rho,y_1,\ldots, y_\tau\}$ such that
\pars
a) \ $v_P^{ }F= v_P^{ }K\oplus\Z v_P^{ }x_1\oplus \ldots\oplus
\Z v_P^{ }x_\rho$,\par
b) \ $y_1P,\ldots,y_\tau P$ is a separating transcendence basis of
$FP|KP$,\par
c) \ $(F,P)$ lies in the absolute inertia field of $(K(T),P)$.
\sn
Assertion c) holds for each transcendence basis $T$ which satisfies
assertions a) and b).
\end{theorem}
\begin{proof}
By Corollary~\ref{fingentb}, the factor group $v_{P}F/v_{P}K$ and the
residue field extension $FP|KP$ are finitely generated. We choose
$x_1,\ldots,x_\rho\in F$ such that $v_P^{ }F=v_{P}K\oplus\Z v_P^{ }
x_1\oplus\ldots\oplus\Z v_P^{ }x_\rho$, where $\rho=\dim_\Q ((v_P^{ }
F/v_P^{ }K)\otimes\Q)$. Since $FP|KP$ is a finitely generated separable
extension, it is separably generated. Therefore, we can choose
$y_1,\ldots,y_\tau \in F$ such that $FP|KP (y_1P, \ldots, y_\tau P)$ is
separable-algebraic, where $\tau=\trdeg FP|KP$. We set $T:=\{x_1,\ldots,
x_\rho,y_1,\ldots,y_\tau\}$ and $F_0:=K(T)$.

 %Inserted by TeXtelmExtel

Now we can choose some $a\in FP$ such that $FP=KP(y_1P,\ldots,y_\tau
P,a)$. Since $a$ is separable-algebraic over $KP(y_1P,\ldots,y_\tau P)$,
by Hensel's Lemma there exists an element $\eta$ in the henselization of
$(F,P)$ such that $\eta P=a$ and that the reduction of the minimal
polynomial of $\eta$ over $F_0$ is the minimal polynomial of $a$ over
$KP(y_1P,\ldots, y_\tau P)$. Then $\eta$ lies in the absolute inertia
field of $F_0\,$. Now the field $F_0(\eta)$ has the same value group and
residue field as $F$, and it is contained in the henselization $F^h$ of
$F\,$. As henselizations are immediate extensions and the henselization
$F_0(\eta)^h$ of $F_0(\eta)$ can be chosen inside of $F^h$, we obtain an
immediate algebraic extension $(F^h| F_0(\eta)^h,P)$. On the other hand,
$(K,P)$ is assumed to be a defectless field. By construction,
$(F_0|K,P)$ is without transcendence defect, and the same is true for
$(F_0(\eta)|K,P)$ since this property is preserved by algebraic
extensions. Hence we know from Theorem~\ref{ai} that $(F_0(\eta),P)$ is
a defectless field. Now Corollary~\ref{corhdie} shows that the extension
$F^h|F_0(\eta)^h$ must be trivial. Therefore, $F$ is contained in
$F_0(\eta)^h$, which in turn is a subfield of the absolute inertia field
of $F_0$. This proves our theorem.
\end{proof}

 %Inserted by TeXtelmExtel

Theorem~\ref{hrwtd} is central in the proof of the uniformization
results presented in this paper. Its importance is based on the fact
that the valuation ring ${\cal O}_{K(T)^i}={\cal O}_P\cap K(T)^i$ is
the strict henselization of the valuation ring ${\cal O}_{K(T)}$ and
therefore the ring extension ${\cal O}_{K(T)^i}|{\cal O}_{K(T)}$ is
local-ind-\'etale; see [Ray], Ch.~X. Using this property, one can
construct an extension $B|A$ of finitely generated $R$-algebras, $R=K$
or $R$ the local base ring appearing in Theorem~\ref{MT2}, such that
$B\subset{\cal O}_P$, $\Frac B =F$, $\Frac A=K(T)$ and $B_q|A_q$ is an
\'etale extension for $q:={\cal M}_{P}\cap A$. In order to prove
Theorems~\ref{MT} and~\ref{MT2}, and ignoring the requirement for the
elements $\zeta\in Z$ for the moment, by the permanence of smoothness
and regularity under \'etale extension it therefore suffices to
construct $A$ in such a way that $A_q$ is smooth over $K$, or is a
regular local ring, respectively. This is done in the next section.

 %Inserted by TeXtelmExtel

%
%---------------------------------------------------------------
%
\section{Abhyankar places on rational function fields} \label{sectrf}
In this section, uniformization of Abhyankar places on rational function
fields of the type appearing in Theorem~\ref{hrwtd} is investigated.
Throughout this section, let $(F|K,P)$ be a rational function field
equipped with a place $P$ subject to the following conditions:

 %Inserted by TeXtelmExtel

\sn
{\it
\textbf{(T)} There exists a transcendence basis $T=(x_1,\ldots,x_\rho,y_1,
\ldots ,y_\tau )$ of $F|K$ such that:
\begin{itemize}
\item $F=K(T)$,
\item $FP=KP(y_1P,\ldots ,y_\tau P)$, with $y_1P,\ldots ,y_\tau P$
algebraically independent over $KP$,
\item $v_P F =v_P K \oplus\Z v_P x_1 \oplus\ldots\oplus\Z v_P x_\rho$.
\end{itemize}
}

 %Inserted by TeXtelmExtel

\noindent
In particular, $P$ is an Abhyankar place of $F|K$.

 %Inserted by TeXtelmExtel

\pars
We fix a finite set $Z\subset{\cal O}_{P}$ and a regular local ring $R
\subseteq K$ such that $\Frac R=K$ and ${\cal M}_{P}\cap R=M$, where $M$
is the maximal ideal of $R$. The case of $R=K$ is included.

 %Inserted by TeXtelmExtel

In the case of $R\ne K$, we cannot prove uniformization over the base
ring $R$ for an arbitrary pair $(P,Z)$ of given data. Instead, we have
to impose rather technical conditions on the pair $(R,Z)$. These
conditions involve the notion of a monoidal transform of $R$: Assume for
the moment that $\dim R>0$ holds, and let $v$ be a valuation of the
field $K$ such that $R\subseteq{\cal O}_{v}$ and $M=R\cap{\cal M}_{v}$.
Let $p\in\Spec R$ be a prime of $\height (p)\geq 1$. The monoidal
transform of $R$ along $v$ with center $p$ is the local ring
\[
R_1:=R[x^{-1}p]_{{\cal M}_{v}\cap R[x^{-1}p]},
\]
where $x\in p$ satisfies $vx=\min \{va \;|\; a\in p\}$. If $p=M$, then
the monoidal transform is also called quadratic transform. It is
well-known that $R_1$ does not depend on the choice of $x$ and that it
is a regular local ring of dimension $\dim R_1\leq\dim R$. Moreover,
$R_1 \subseteq {\cal O}_{v}$ and $M_1=R_1 \cap{\cal M}_{v}$ for the
maximal ideal $M_1$ of $R_1$. Every member of a finite chain $R=:R_0\leq
R_1\leq \ldots\leq R_t$ of local rings, where $R_{i+1}$ is a monoidal
transform of $R_i$ along $v$ with center $p_i\in\Spec R_i$, is called an
iterated monoidal transform of $R$ along $v$.

 %Inserted by TeXtelmExtel

\pars
The properties of $(R,Z)$ we are interested
in can now be formulated as follows:
\sn
{\it
\textbf{(NC)}
There exists an iterated monoidal transform $R^\prime$ of $R$ along the
valuation $v_{P}|_{K}$ such that the regular local ring $R^\prime$
admits a regular parameter system $(t_1,\ldots ,t_d)$ with the following
property: every $\zeta\in Z$ admits a representation
\begin{equation}                              \label{first fraction}
\zeta \>=\>\frac{\sum\limits_{i=1}^N
r_i\,\underline{x}^{\underline{\mu}_i}\underline{y}^{\underline{\nu}_i}}
{\sum\limits_{i=1}^N s_i\,\underline{x}^{\underline{\kappa}_i}
\underline{y}^{\underline{\lambda}_i}}\>,\;\; r_i,s_i\in R^\prime,
\end{equation}
with $\underline{x}:=(x_1,\ldots ,x_\rho)$, $\underline{y}:=(y_1,\ldots
, y_\tau)$, $\underline{\mu}_i$, $\underline{\kappa}_i\in\N_0 ^\rho$,
$\underline{\nu}_i$, $\underline{\lambda}_i\in\N_0^\tau$, where the
prime factorizations in $R^\prime$ of the coefficients $r_i, s_i$ have
the form
\begin{equation}                              \label{normal crossing}
u\prod\limits_{i=1}^d t_i^{\epsilon_i}, \;\; u\in (R^\prime ) ^\times ,
\;\epsilon_i\in\N_0\> ,
\end{equation}
i.e., the coefficients $r_i,s_i$ are $R^\prime$-momomials in $\{ t_1,\ldots ,t_d\}$.
}

 %Inserted by TeXtelmExtel

\sn
{\it
\textbf{(V)}
The values $v_P t_1,\ldots ,v_P t_\delta$ of those parameters $t_1,\ldots
,t_\delta$ actually occurring in at least one of the prime factorizations
(\ref{normal crossing}), are rationally independent. Without loss of
generality, we assume here that these parameters are the first $\delta$ of
the complete set.
}

 %Inserted by TeXtelmExtel

\pars
The main result of this subsection now reads as follows:
\begin{theorem}                            \label{uratff}
Let $(F|K,P)$ be a valued rational function field satisfying the
requirement (T). Let $R\subseteq K$ be a universally catenary, regular
local ring dominated by ${\cal O}_{P}$ and let $Z\subset{\cal O}_{P}$ be
a finite set such that the pair $(R,Z)$ fulfills the requirements (NC)
and (V). Then there exist an iterated monoidal transform $R^\prime$ of
$R$ along $v_{P}|_{K}$ and elements $x_1^\prime,\ldots
,x_{\rho+\delta}^\prime\in{\cal O}_{P}$, $\delta\leq\dim R^\prime$, such
that the localization of the $R^\prime$-algebra $A:=R^\prime
[x_1^\prime,\ldots ,x_{\rho+\delta}^\prime,y_1,\ldots ,y_\tau]$ at the
prime $q:={\cal M}_{P} \cap A$ is a regular ring having the properties:
$\dim A_q= \dim R^\prime +\rho$, the elements $x_1^\prime,\ldots
,x_{\rho+\delta}^\prime$ are a part of a regular parameter system of
$A_q$, $Z\subset A_q$ and every element of $Z$ is an $A_q$-monomial in
$\{x_1^\prime,\ldots ,x_{\rho+\delta}^\prime\}$.

 %Inserted by TeXtelmExtel

In particular, there exists an $R$-model $X$ of $F$ such that $P$ is
centered in a regular point $x\in X$ with the properties $Z\subset
{\cal O}_{X,x}$ and $\dim {\cal O}_{X,x}\leq\dim R+\rho$. Moreover, if
$R=K$ then the model $X$ can be chosen such that $X\cong\A_K^{\rho
+\tau}$ and $\dim {\cal O}_{X,x}=\rho$.
\end{theorem}

 %Inserted by TeXtelmExtel

The following lemma is applied in the proof of Theorem~\ref{uratff}; it
was proved (but not explicitly stated) by Zariski in [Z] for subgroups
of $\R$, using the algorithm of Perron. We leave it as an easy exercise
to the reader to prove the general case by induction on the rank of the
ordered abelian group. However, an instant proof of the lemma can also
be found in [El] (Theorem~2.2).
\begin{lemma}                               \label{perron}
Let $\Gamma$ be a finitely generated ordered abelian group. Take any
non-negative elements $\alpha_1,\ldots,\alpha_\ell\in\Gamma$. Then there
exist positive elements $\gamma_1,\ldots,\gamma_{\rho}\in\Gamma$ such
that $\Gamma=\Z \gamma_1 \oplus\ldots\oplus\Z\gamma_{\rho}$ and every
$\alpha_i$ can be written as a sum $\sum_{j}^{} n_{ij}\gamma_j$ with
non-negative integers $n_{ij}\,$.
\end{lemma}

 %Inserted by TeXtelmExtel

We now turn to the
\sn
{\bf Proof of Theorem~\ref{uratff}}: \
Since iterated monoidal transforms of a local ring $R$ are essentially of
finite type over $R$, the existence of the scheme $X$ is a consequence
of the first part of the theorem.

 %Inserted by TeXtelmExtel

To simplify notation, we replace $R$ with an iterated monoidal transform
$R^\prime$ of $R$ along the valuation $v_{P}|_{K}$ having the properties
(NC) and (V) for a specific regular parameter system $t_1,\ldots ,t_d\in
M$, $M$ the maximal ideal of $R^\prime$. Recall that along with $R$ the
domain $R^\prime$ is universally catenary and that $\dim
R^\prime\leq\dim R$.

 %Inserted by TeXtelmExtel

The coefficients appearing in the representations (\ref{first fraction})
can now be replaced by their prime factorizations in $R$:
\begin{equation}                              \label{second fraction}
\zeta \>=\>\frac{\sum\limits_{i=1}^N
u_i\,\underline{t}^{\underline{\epsilon}_i}
\underline{x}^{\underline{\mu}_i}\underline{y}^{\underline{\nu}_i}}
{\sum\limits_{i=1}^N v_i \,\underline{t}^{\underline{\delta}_i}
\underline{x}^{\underline{\kappa}_i}
\underline{y}^{\underline{\lambda}_i}}\>,\;\; u_i,v_i\in R^\times,
\end{equation}
where $\underline{t}=(t_1,\ldots ,t_\delta)$ are those regular parameters
among the $t_1,\ldots ,t_d$ that actually occur in at least one of the
prime factorizations (\ref{normal crossing}) of the coefficients.
Note that in our notation in (\ref{second fraction}) the dependence of the
coefficients and exponents on $\zeta$ does not explicitly appear, in
order to avoid overloading. We use this simplification throughout the
present proof.

 %Inserted by TeXtelmExtel

Next, one divides numerator and denominator by a monomial in
$\underline{t}$ and $\underline{x}$ with least $v_P$-value among the
monomials appearing in the denominator. One can assume that this
monomial is the first one in the denominator, thus obtaining the
expression:
\begin{equation}
\label{third fraction}
\zeta \>=\>\frac{\sum\limits_{i=1}^N
u_i\,\underline{t}^{\underline{\epsilon}_i-
\underline{\delta}_1}\underline{x}^{\underline{\mu}_i-\underline{\kappa}_1}
\underline{y}^{\underline{\nu}_i}}{\sum\limits_{i=1}^N v_i\,
\underline{t}^{\underline{\delta}_i
-\underline{\delta}_1}\underline{x}^{\underline{\kappa}_i-
\underline{\kappa}_1}\underline{y}^{\underline{\lambda}_i}}\>,\;\;
u_i,v_i\in R^\times .
\end{equation}
By construction, the monomial expressions
$\underline{t}^{\underline{\delta}_i
-\underline{\delta}_1}\underline{x}^{\underline{\kappa}_i
-\underline{\kappa}_1}$ all have non-negative $v_P$-value. Since
$v_P(\zeta )\geq 0$ and since $P$ is an Abhyankar place the same holds
for the monomial expressions $\underline{t}^{\underline{\epsilon}_i
-\underline{\delta}_1}\underline{x}^{\underline{\mu}_i
-\underline{\kappa}_1}$ appearing in the numerator; see
Theorem~\ref{prelBour}.

 %Inserted by TeXtelmExtel

Among the monomials in $\underline{t}$ and $\underline{x}$ appearing
in the numerator of (\ref{third fraction}) we choose the unique one having
least $v_P$-value, $h_{\zeta}:=\underline{t}^{\underline{\epsilon}_m-\underline{\delta}_1}
\underline{x}^{\underline{\mu}_m-\underline{\kappa}_1}$ say, and define
the finite sets
\[
\begin{array}{rcl}
H_{\zeta}\;:=\;\{h_{\zeta}\}&\cup&
\{\underline{t}^{\underline{\delta}_i -\underline{\delta}_1}
\underline{x}^{\underline{\kappa}_i-\underline{\kappa}_1}
\; |\; i=1,\ldots ,N\}\\
&\cup&
\{\underline{t}^{\underline{\epsilon}_i-\underline{\epsilon}_m}
\underline{x}^{\underline{\mu}_i-\underline{\mu}_m}
=\frac{\underline{t}^{\underline{\epsilon}_i-\underline{\delta}_1}
\underline{x}^{\underline{\mu}_i-\underline{\kappa}_1}}
{h_{\zeta}}\mid i=1,\ldots ,N\}.
\end{array}
\]
The finite set
\[
H:=\{t_1,\ldots ,t_\delta,x_1,\ldots ,x_\rho\}\cup\bigcup
\limits_{\zeta\in Z}H_{\zeta}
\]
then consists of elements with  non-negative $v_P$-value only.

 %Inserted by TeXtelmExtel

Let $G\subset K(x_1,\ldots ,x_\rho )^\times$ be the group (freely)
generated by the elements $t_1,\ldots ,t_\delta$ and $x_1,\ldots
,x_\rho$. By assumption, the valuation $v_P$ induces an isomorphism
\[
v_P:\; G\>\rightarrow\>\bigoplus\limits_{i=1}^\delta\Z v_P t_i\oplus
\bigoplus \limits_{i=1}^\rho\Z v_P x_i\subseteq v_{P}F.
\]
Applying Lemma~\ref{perron} yields a basis $(v_P x_1^\prime ,\ldots ,
v_P x_{\rho+\delta}^\prime )$ of the group $v_PG$ consisting of positive
elements such that every $v_P h \in v_P H\subset v_P G$ can be expressed as
a linear combination with non-negative coefficients. It follows that
every $h\in H$ can be expressed in the form
\begin{equation}
\label{monomials}
h\>=\>\prod\limits_{i=1}^{\rho +\delta}{x_i^{\prime}}^{\>\mu_i^\prime},\;
\mu_i^\prime\in\N_0\>.
\end{equation}
By definition of the set $H_{\zeta}$ the monomials in $\underline{t}$ and
$\underline{x}$ appearing in the numerator of (\ref{third fraction}) have
the form $hh_{\zeta}$ for some $h\in H_{\zeta}$, therefore every of these
monomials has the form
\begin{equation}
\label{monomial product}
\underline{t}^{\underline{\epsilon}_i-\underline{\delta}_1}
\underline{x}^{\underline{\mu}_i-\underline{\kappa}_1}
\>=\>(\prod\limits_{i=1}^{\rho +\delta}{x_i^{\prime}}^{\>\beta_i})
(\prod\limits_{i=1}^{\rho +\delta}{x_i^{\prime}}^{\>\alpha_i}),\;
\beta_i\in\N_0\>,
\end{equation}
where
\begin{equation}
\label{minimal monomial}
h_{\zeta}\>=\>\prod\limits_{i=1}^{\rho +\delta}
{x_i^{\prime}}^{\>\alpha_i},\;\alpha_i\in\N_0\>.
\end{equation}
Substituting (\ref{monomials}) and (\ref{monomial product}) in (\ref{third fraction}) yields
\begin{equation}
\label{fourth fraction}
\zeta \>=\>\frac{\sum\limits_{i=1}^N u_i\,\underline{x^\prime}^{\>\underline{\beta}_i}\,
\underline{y}^{\underline{\nu}_i}}
{\sum\limits_{i=1}^N v_i\,\underline{x^\prime}^{\>\underline{\kappa}_i^\prime}\,
\underline{y}^{\underline{\lambda}_i}}
\;(\prod\limits_{i=1}^{\rho +\delta}
{x_i^{\prime}}^{\>\alpha_i})
\>,\;\; u_i,v_i\in R^\times ,
\end{equation}
with non-negative exponents $\underline{\beta}_i ,
\underline{\kappa}_i^\prime\in\N_0^{\rho +r}\,$.

 %Inserted by TeXtelmExtel

We set $A:=R[x_1^\prime,\ldots ,x_{\rho+\delta}^\prime,y_1,\ldots ,y_\tau ]
\subseteq{\cal O}_{P}$ and claim: $A_{q_A}$, where $q_A:={\cal M}_{P}
\cap A$, is a regular ring of dimension $\rho +\dim R$, $Z\subset A_{q_A}$
and all $\zeta\in Z$ are $A_{q_A}$-monomials in
$\{x_1^\prime,\ldots ,x_{\rho+\delta}^\prime\}$.

 %Inserted by TeXtelmExtel

Let $B:=R[x_1^\prime,\ldots, x_{\rho+\delta}^\prime]\subseteq A$ and
consider the ideals
\[
J:=\sum\limits_{i=1}^{\rho +\delta}Bx_i^\prime
+\sum\limits_{i=\delta +1}^d Bt_i
\subseteq q_B:={\cal M}_{P}\cap B.
\]
Due to the equations (\ref{monomials}) for the elements $t_i$,
$i=1,\ldots ,\delta$, one gets $J\cap R =M$. Hence, $B/J=R/M$ and
therefore, $q_B=J$ is a maximal ideal of $B$ generated by $\rho+d$
elements, $d=\dim R$.

 %Inserted by TeXtelmExtel

According to the altitude formula in finitely generated integral domains
over universally catenary rings, one has:
\[
\begin{array}{rl}
\height q_B & =\height  M +\trdeg (B|R) -\trdeg (B/q_B |R/M)\\
&=d+\rho,
\end{array}
\]
where we use that due to the equations (\ref{monomials}) for the
elements $x_i$ the relation $\Frac B=K(x_1,\ldots ,x_\rho)$ holds.
Consequently, $B_{q_B}$ is a regular local ring of Krull dimension
$d+\rho$.

 %Inserted by TeXtelmExtel

Considering that $A_{q_A}=B_{q_B}[y_1,\ldots ,y_\tau ]_p$, where
$p={\cal M}_{P} \cap B_{q_B}[y_1,\ldots ,y_\tau ]$, and that
$B_{q_B}[y_1,\ldots ,y_\tau ]$ is a polynomial ring over $B_{q_B}$, we
get the desired regularity of $A_{q_A}$ because a polynomial ring over a
regular ring is regular.

 %Inserted by TeXtelmExtel

It remains to calculate the dimension of $A_{q_A}$. By the choice of the
transcendence basis $T$ we have that  $B_{q_B}[y_1,\ldots ,y_\tau ]/p=
R/M [y_1+p,\ldots ,y_\tau +p]$ is a polynomial ring in the variables
$y_i+p$, $i=1,\ldots ,\tau$. Thus $p=q_B[y_1,\ldots ,y_\tau]$ and
therefore, $\height p=\height q_B$, where the latter equation holds in
any polynomial ring over a noetherian ring. Consequently, $A_{q_A}$ is a
local ring of Krull dimension $\rho+d$.

 %Inserted by TeXtelmExtel

If $R=K$ is a field, then the $K$-algebra $B$ is a polynomial
ring in the variables $x_1^\prime ,\ldots ,x_\rho^\prime$ and therefore,
$A_{q_A}$ is a localization of the polynomial ring $K[x_1^\prime ,\ldots
,x_\rho^\prime,y_1, \ldots ,y_\tau]$.

 %Inserted by TeXtelmExtel

Turning to the claim about $Z$ we first observe that the denominator of the
expression (\ref{fourth fraction}) has $v_P$-value $0$, which implies the
inclusion $Z\subset A_{q_A}$. The numerator of the first factor in
(\ref{fourth fraction}) has $v_P$-value $0$ too, because its $m$th summand
has $v_P$-value $0$. Thus this first factor is a unit in $A_{q_A}$, which verifies
the claim.
\QED

 %Inserted by TeXtelmExtel

 %Inserted by TeXtelmExtel

We finish this section with a discussion of the requirements (NC) and (V)
in various cases to give examples in which Theorem~\ref{uratff} applies.

 %Inserted by TeXtelmExtel

\sn
\textbf{(1)} $\dim R\leq 1$: In this case, $R$ is either the field $K$
itself or a discrete valuation ring. In the first case, the coefficients
appearing in the representations (\ref{first fraction}) are units, hence
(NC) and (V) are trivially satisfied. In the second case, every prime
element $t\in R$ satisfies the requirements (NC) and (V), so
$R^\prime=R$.

 %Inserted by TeXtelmExtel

\sn
\textbf{(2)} ${\cal O}_{K}$ is a discrete valuation ring: It is known
that in this case, ${\cal O}_{K}=\cup_{k=0}^\infty R_k$, where $R_k$ is
the quadratic transform of $R_{k-1}$ along $v_{P}|_{K}$ --- see for
example [HRW]. Given a finite set $Z^\prime\subset{\cal O}_{K}\,$, one
chooses prime factorizations $z=u_zt^{\epsilon_z}$, $u_z\in {\cal
O}_{K}^\times$, $t$ a prime element, for all $z\in Z^\prime$. It is now
easy to verify that these prime factorizations remain valid in some
$R_\ell\,$, thus showing that (NC) and (V) are satisfied
for arbitrary finite sets $Z\subset {\cal O}_{P}$.

 %Inserted by TeXtelmExtel

\sn
\textbf{(3)} $\dim R=2$: One applies Abhyankars results [A6]. If $\trdeg
(KP|R/M)=1$, then the valuation $v_{P}|_{K}$ is discrete, thus we are in
case (2). In the case of $\trdeg (KP|R/M)=0$, due to [A6], \S 4, Thm.~2,
for every finite set $Z^\prime\subset {\cal O}_{K}$ there exists a
$2$-dimensional monoidal transform $R^\prime$ with a regular parameter
system $(t_1,t_2)$ such that the elements $z\in Z^\prime$ have prime
factorizations of the type (\ref{normal crossing}) with property (V),
where $\delta\in\{1,2\}$ depending on the rational rank of $v_{P}K$. In
the case of $\dim {\cal O}_{K}=2$ one has to assume in addition that $R$
is a Nagata ring. It follows that Theorem~\ref{uratff} applies for
$2$-dimensional regular local Nagata rings $R$.

 %Inserted by TeXtelmExtel

\sn
\textbf{(4)} $\dim R=3$: In this case, the following two results are
relevant in the current context:

 %Inserted by TeXtelmExtel

If $R$ is excellent and $\chara R=\chara R/M$ holds, then every finite
set $Z^\prime\subset{\cal O}_{K}$ is contained in an iterated monoidal
transform $R^\prime$ of $R$ along $v_{P}|_{K}$ such that the elements
$z\in Z^\prime$ have prime factorizations of the form (\ref{normal
crossing}) - [A3], (5.2.3). Thus (NC) is satisfied for every finite
set $Z\subset{\cal O}_{P}$.

 %Inserted by TeXtelmExtel

D. Fu proved ([F], Prop. 3.5) that (NC) and (V) are satisfied for
every finite set $Z\subset{\cal O}_{P}\,$, provided that $R$ is
essentially of finite type over an algebraically closed field $k$ and
$KP=k$.

 %Inserted by TeXtelmExtel

%
%---------------------------------------------------------------
%
\section{The general case}                        \label{general}
In this section we provide the proofs for Theorems~\ref{MT}
and~\ref{MT2}. The ingredients are already given in Sections~\ref{sect3}
and~\ref{sectrf}. We need one more fact --- an ascend property for
$R$-uniformizability --- to bring them together.

 %Inserted by TeXtelmExtel

Consider a finite extension $(F|E,P)$ of valued fields such that $F$ is
contained in the absolute inertia field of $(E,P|_{E})$. It is then
well-known that the extension ${\cal O}_{F}|{\cal O}_{E}$ is
local-\'etale ([Ray], Ch.~X, Thm.~1): ${\cal O}_{F}=A_q$ for an \'etale
${\cal O}_{E}$-algebra $A$ and $q=A\cap {\cal M}_F$. According to [Ray],
Ch.~V, Thm.~1 we can assume that $A$ is standard-\'etale, i.e.,
\begin{equation}
\label{standard}
A={\cal O}_{E}[x]_{g(x)},
\end{equation}
where ${\cal O}_{E}[x]={\cal O}_{E}[X]/f{\cal O}_{E}[X]$ with a monic
polynomial $f\in{\cal O}_{E}[X]$. Furthermore, $g\in{\cal O}_{E}[X]$ is
chosen such that the image of the derivative $f^\prime$ under the
natural morphism $\phi:{\cal O}_{E}[X]\rightarrow A$ is a unit.

 %Inserted by TeXtelmExtel
\pars
\textbf{Claim:} In the definition of $A$ we can assume $f$ to be
prime.

 %Inserted by TeXtelmExtel

\noindent
The Lemma of Gau{\ss} allows to factorize $f$ as $\prod\limits_{i=1}^r
p_i^{\nu_i}$ with pairwise distinct, monic prime polynomials
$p_i\in{\cal O}_{E}[X]$. Among them there is a unique $p_j$ with $\phi
(p_j)=p_j(x)=0$. We consider the natural surjection $\psi :A\rightarrow
({\cal O}_{E}[X]/p_j{\cal O}_{E}[X])_{g(x)}$. The equation $\psi
(\frac{h+f{\cal O}_{E}[X]}
{g(x)^s})=0$ implies $h\in p_j{\cal O}_{E}[X]$. Since in $A$ we have
$p_j(x)=0$ there exists $t\in\N$ such that $g^t p_j\in f{\cal O}_{E}[X]$,
thus $g^t h\in f{\cal O}_{E}[X]$ and hence $\frac{h+f{\cal O}_{E}[X]}
{g(x)^s}=0$. Finally, for some $p^\star\in{\cal O}_{E}[X]$ we have $\phi
(f^\prime )=\phi ((p^\star p_j)^\prime )=\phi (p^\star p_j ^\prime )\in
A^\times$, which implies that $\phi (p_j^\prime )\in A^\times$. This
proves that $\psi$ is an isomorphism and hence the claim.

 %Inserted by TeXtelmExtel

Next, we fix a set of structural constants determining $A$ uniquely,
which we shall use to define an \'etale algebra over a subring of
${\cal O}_{E}\,$. Let $h\in {\cal O}_{E}[X]$ be chosen such that
$f^\prime(x)^{-1}=\frac{h(x)} {g(x)^s}$, $s\in\N$, and let
$C(f,g,h)\subset{\cal O}_{E}$ be the set of coefficients of the
polynomials $f, g$ and $h$.

 %Inserted by TeXtelmExtel

Let $Z\subset {\cal O}_F$ be a finite set and split it as
$Z=(Z\cap {\cal O}_F^\times )\cup  (Z\cap {\cal M}_F )$. Since $v_P$ is
unramified in the extension $F|E$ we can write every $\zeta\in Z\cap
{\cal M}_F$ in the form $\zeta =u_{\zeta}\zeta^\prime$ with $u_{\zeta}\in
{\cal O}_F^\times$ and $\zeta^\prime\in {\cal M}_E$; let $Z^\prime
:=\{\zeta^\prime\;|\;\zeta\in Z\cap {\cal M}_F\}$.

 %Inserted by TeXtelmExtel

For each $\xi\in Z^\times:=(Z\cap {\cal O}_F^\times )\cup\{u_{\zeta}\;|\;
\zeta\in Z\cap {\cal M}_F\}$ we choose a representation
\begin{equation}
\label{repres}
\xi =\frac{a(x)}{b(x)}g(x)^k,\;\; a,b\in{\cal O}_E[X], \; b(x)\not\in
q,\; k\in\Z.
\end{equation}
Finally we define the finite set  $C(Z^\times )\subset {\cal O}_E$ as the
collection of the coefficients of the polynomials $a,b$ appearing in the
representations (\ref{repres}) for the elements $\xi$.
\begin{lemma}                                   \label{etale}
Let $(F|E,P)$ be a finite extension of valued fields such that $F$ is
contained in the absolute inertia field $E^i$ of $(E, P|_{E})$. Let
$Z\subset{\cal O}_{F}$ be a finite set and let $R\subseteq {\cal O}_{E}$
be a regular local ring with maximal ideal $M={\cal M}_{E}\cap R$.

 %Inserted by TeXtelmExtel

Let $A\subseteq{\cal O}_{F}$ be the \'etale ${\cal O}_E$-algebra
(\ref{standard}) and assume that there exists a set of re\-presentations
(\ref{repres}) for the elements $\xi\in Z^\times$, such that the pair
$(P|_{E},C(f,g,h)\cup C(Z^\times )\cup Z^\prime )$ is $R$-uniformizable.
Then the pair $(P,Z)$ is $R$-uniformizable, too.

 %Inserted by TeXtelmExtel

Moreover, one can find $R$-models $X$ of $E$ and $Y$ of $F$ and a
morphism $\pi :\; U\rightarrow V$, where $U$ is an affine open neighborhood
of the center $y$ of $P$ on $Y$, and $V$ is an affine open neighborhood
of the center $x$ of $P|_E$ on $X$ such that: $x,y$ are regular points
on the respective model, $\pi y=x$, and the extension ${\cal O}_{Y,y}|
{\cal O}_{X,x}$ is local-\'etale.

 %Inserted by TeXtelmExtel

In particular, $\dim {\cal O}_{Y,y}=\dim {\cal O}_{X,x}$ holds.
Furthermore, if $R=K$ is a field and ${\cal O}_{X,x}$ is smooth over
$K$, then ${\cal O}_{Y,y}$ is smooth over $K$.

 %Inserted by TeXtelmExtel

The $R$-model $X$ can be chosen such that for some regular system of
parameters\linebreak
$(a_1,\ldots ,a_d)$ of ${\cal O}_{X,x}$ every $\zeta^\prime\in Z^\prime$
is an ${\cal O}_{X,x}$-monomial in $\{ a_1,\ldots ,a_d\}$. In this case
every $\zeta\in Z$ is an ${\cal O}_{Y,y}$-monomial in the regular system
of parameters $(a_1,\ldots ,a_d)$ of ${\cal O}_{Y,y}$.
\end{lemma}
\begin{proof}
By assumption there exists a finitely generated $R$-algebra
$B\subseteq{\cal O}_E$ such that $B_{q_B}$, $q_B:={\cal M}_E\cap B$, is
regular and contains the finite set $C(f,g,h)\cup C(Z^\times )\cup Z^\prime$.
Define $C:=B[x]_{g(x)}\subseteq A$, where $x$ is the element appearing in
the definition (\ref{standard}) of $A$. We have $\Frac C=F$, and the
$B_{q_B}$-algebra $C_{q_B}=B_{q_B}[x]_{g(x)}$ is standard-\'etale: this is
a consequence of the construction of $C_{q_B}$ once we have verified
that $B_{q_B}[x]\cong B_{q_B}[X]/f B_{q_B}[X]$. So assume $h(x)=0$ for
some $h\in B_{q_B}[X]$; since $f$ is the minimal polynomial of $x$ over
$E$ we get $h=fh^\star$, $h^\star\in E[X]$. Now $B_{q_B}$ is integrally
closed in $E$ and $f$ is monic, thus the Lemma of Gau{\ss} yields
$h^\star\in B_{q_B}[X]$.

 %Inserted by TeXtelmExtel

Since regularity ascends in \'etale extensions, the domain $C_{q_B}$ and
hence also $C_{q_C}$, where $q_C:={\cal M}_F\cap C$, are regular.
%Since regularity ascends in \'etale extensions the domain $C_{q_B}$ and
%in particular $C_{q_C}$, $q_C:={\cal M}_F\cap C$, is regular.
Moreover, we have $Z\subset C_{q_C}$ by construction. \'Etale
extension preserves the Krull dimension and smoothness. Moreover, if
$(a_1,\ldots ,a_d)$ is a regular parameter system of $B_{q_B}\,$, then
it is a regular parameter system for $C_{q_C}$, too. These facts yield
the remaining assertions.
\end{proof}

 %Inserted by TeXtelmExtel

We are now prepared to prove our main results.

 %Inserted by TeXtelmExtel

\mn
\textbf{Proof of Theorems~\ref{MT} and~\ref{MT2}:}
\sn
One starts by choosing a transcendence basis $T\subset F$ with the
properties described in Theorem~\ref{hrwtd}; in particular, the valued
field $(F,P)$ lies in the absolute inertia field of $(K(T),P|_{K(T)})$.

 %Inserted by TeXtelmExtel

According to Lemma~\ref{etale}, the pair $(P,Z)$ is $R$-uniformizable
for a given finite set $Z\subset {\cal O}_{P}$ once the pair
$(P|_{K(T)}, Z^{\prime\prime})$ is $R$-uniformizable for a certain
finite set $Z^{\prime\prime}\subset {\cal O}_{K(T)}$ derived from $Z$
and the extension ${\cal O}_{P}|{\cal O}_{K(T)}$ --- see the discussion
preceding Lemma~\ref{etale}. Moreover, the elements $\zeta\in Z$ possess
the required factorization property once the elements of a certain
subset $Z^\prime\subseteq Z^{\prime\prime}$ possess this factorization
property.

 %Inserted by TeXtelmExtel

The valued rational function field $(K(T),P|_{K(T)})$ satisfies the
requirements (T) of Section~\ref{sectrf}. Points (1) and (3) of the
discussion at the end of Section~\ref{sectrf} show that the pair
$(R,Z^{\prime\prime})$ fulfills the requirements (NC) and (V) for {\em every}
finite set $Z^{\prime\prime}\subset {\cal O}_{K(T)}$. An application of
Theorem~\ref{uratff} thus yields $R$-uniformizability of $(P|_{K(T)},
Z^{\prime\prime} )$ and the factorization property for the elements $\zeta
\in Z^{\prime\prime}$ for an arbitrary finite set $Z^{\prime\prime}$.

 %Inserted by TeXtelmExtel

In the case of $R=K$, the pair $(P|_{K(T)},Z^{\prime\prime} )$ is
$K$-uniformizable on the affine space $X=\A_K^{\rho +\tau}$, $\rho
=\dim_{\Q}(v_P F\otimes\Q )$, due to Theorem~\ref{uratff}. Moreover, the
center $x\in X$ of $P|_{K(T)}$ satisfies $\dim {\cal O}_{X,x}=\rho$.
Lemma~\ref{etale} thus gives the remaining assertions of
Theorem~\ref{MT}.

 %Inserted by TeXtelmExtel

In the case of $R\ne K$, it remains to verify the dimension statements
of Theorem~\ref{MT2}. If $R$ or ${\cal O}_K$ is a discrete valuation
ring, then points (1) and (2) at the end of Section~\ref{sectrf}
combined with Lemma~\ref{etale} give the assertion. In the remaining
case, point (3) of Section~\ref{sectrf} contains the relevant
information.
\QED

 %Inserted by TeXtelmExtel

As an immediate corollary of Theorem~\ref{MT} we get:
\begin{corollary}                                  \label{Cor MT}
Assume the situation of Theorem~\ref{MT} to be present except for the
separability of  $FP|K$. Then there exists a finite purely inseparable
extension $L|K$ such that $(\widehat{P},Z)$, $\widehat{P}$ the unique
prolongation of $P$ to the constant extension $FL|L$, is
$L$-uniformizable. All other assertions of Theorem~\ref{MT} remain valid
over $L$.
\end{corollary}

 %Inserted by TeXtelmExtel

Instead of stating the mere existence of an algebraic or arithmetic
variety $X$ on which the given place $P$ can be uniformized as in
Theorems~\ref{MT} and~\ref{MT2}, one can rather explicitly describe the
structure of an affine scheme $X$ which does the job. The description
follows directly from the proof of the two theorems. In the case of a
$K$-trivial place $P$, one can view this description as a structure
theorem for the valued function field $(F|K,P)$.
\begin{theorem}                                      \label{structMT}
Let $F|K$ be a function field and $P$ an Abhyankar place of $F|K$ such
that $(K,P)$ is defectless, $FP|KP$ is separable and the group
$v_{P}F/v_{P}K$ is torsion-free. Further, let $R\subset K\cap {\cal O}_{P}$, $\Frac R=K$, be a
noetherian, regular local ring with $\dim R\leq 2$ and maximal ideal
$M={\cal M}_{P}\cap R$. Assume that $R$ is a Nagata ring if $\dim R=2$.
Let $Z\subset {\cal O}_{P}$ be a finite set.

 %Inserted by TeXtelmExtel

Then there exists a transcendence basis $T^\prime=(x_1^\prime,\ldots ,x_\rho^\prime
,y_1,\ldots , y_\tau)\subset {\cal O}_{P}$ of $F|K$, an iterated
monoidal transform $R^\prime$ of $R$ along $v_{P}|_{K}$ and a finitely
generated $R^\prime [T^\prime]$-algebra $A\subset {\cal O}_{P}$ having the
properties:
\begin{enumerate}
\item $y_1P,\ldots ,y_\tau P$ form a separating transcendence basis
of $FP|KP$.
\item $v_{P}x_1^\prime,\ldots ,v_{P}x_\rho^\prime$ are rationally independent elements.
\item $\Frac A=F$ and $Z\subset A_q$, where $q:={\cal M}_{P}\cap A$.
\item $A_q$ is a regular local ring of Krull dimension
$\rho +\dim R^\prime $.
\end{enumerate}
If $P$ is trivial on $K$ the elements appearing in point 2 form a basis of $v_{P}F$ over $\Z$.
Moreover in this case the extension $A_q|K[T^\prime ]_{q\cap K[T^\prime ]}$ is local-\'etale.
\end{theorem}

 %Inserted by TeXtelmExtel

%
%---------------------------------------------------------------
%
\section{Uniformization of finitely many Abhyankar places on an affine
scheme}
The main result of the present section is a proof of Theorem~\ref{MT3}.
We start with some valuation theory that eventually allows to take over
the proof of Theorems~\ref{MT} and~\ref{MT2}.

 %Inserted by TeXtelmExtel

Throughout this section, we assume that the hypotheses of
Theorem~\ref{MT3} hold. We set $v_j= v_{P_j}\,$, $j=1,\ldots ,\ell$, and
denote the valuation ring of $v_j$ by ${\cal O}_{j}$. Moreover, we set
${\cal O}:=\cap_{j=1}^\ell{\cal O}_{j}$. We may assume that the valuation
rings ${\cal O}_{j}$ are pairwise incomparable with respect to inclusion.
This follows from a more general observation:

 %Inserted by TeXtelmExtel

\begin{lemma}                              \label{coars}
Let $(F|K,P)$ be a valued function field and $R\subseteq{\cal O}_{K}$ a
regular local ring. Let $Q$ be a \bfind{coarsening of} $P$, i.e., a
place of $F$ whose valuation ring contains that of $P$. If the pair
$(P,Z)$, $Z\subset {\cal O}_{P}$ finite, is $R$-uniformizable on
the scheme $X$, then the same is true for the pair $(Q,Z)$.
\end{lemma}
\begin{proof}
By assumption, $P$ is centered in a regular point $x$ of an $R$-scheme
$X$ of finite type and $Z\subset {\cal O}_{X,x}\subseteq {\cal O}_{P}$.
It follows that for $q:={\cal M}_{Q}\cap {\cal O}_{X,x}$ the local ring
$({\cal O}_{X,x})_q$ is regular and contains $Z$.
\end{proof}
\n
Consider now the set $\{P_1,\ldots ,P_{\ell}\}$ and a finite set $Z\subset
\cap_{j=1}^\ell {\cal O}_{j}$. Let $\{P_1,\ldots ,P_{\ell^\prime}\}$ be the
subset of those places $P_j$ for which ${\cal O}_{j}$ is minimal in the
set $\{{\cal O}_{1}, \ldots , {\cal O}_{\ell}\}$, with respect to
inclusion. Assume that Theorem~\ref{MT3} is proved for this smaller set
of places and the finite set $Z$. Let $x_1,\ldots ,x_{\ell^\prime}$ be
the regular centers of the places $P_j$, $j=1,\ldots , \ell^\prime$, on
the scheme $X$. Then applying Lemma~\ref{coars} we see that the places
$P_j$, $j=\ell^\prime+1,\ldots ,\ell$, are centered in regular points
$x_j\in X$, $j=\ell^\prime+1,\ldots ,\ell$, with the property that their
local rings are localizations of some of the local rings
${\cal O}_{X,x_j}$, $j=1,\ldots ,\ell^\prime$. Therefore we have
$Z\subset\cap_{j=1}^\ell {\cal O}_{X,x_j}$ as required in
Theorem~\ref{MT3}. {\em Thus, from now on we assume the rings
${\cal O}_{j}$ to be pairwise incomparable.}

 %Inserted by TeXtelmExtel

The value groups $v_j F$ are assumed to be isomorphic over $v_P^{ } K$.
Therefore, passing to equivalent valuations if necessary, we may assume
that $v_1 F=\ldots=v_\ell F$. Since this group is finitely generated over
$vK$ by Corollary~\ref{fingentb}, we can write $v_1F=v_P^{ }
K\oplus\bigoplus_{i=1} ^{\rho}\Z \alpha_i\,$ with positive $\alpha_i\,$.

 %Inserted by TeXtelmExtel

By assumption, we have $\tau =\trdeg (FP_j |KP)$, $j=1,\ldots ,\ell$. If
$\tau >0$, then we choose a separating transcendence basis $\xi_1,
\ldots, \xi_\tau$ of the separable extension $FP_1|KP$. By passing to
equivalent places if necessary, we may assume that $\xi_1,\ldots,
\xi_\tau$ is also a separating transcendence basis of the separable
extensions $FP_j|KP$, $j=2,\ldots ,\ell$.

 %Inserted by TeXtelmExtel

\pars
In this situation we prove:

 %Inserted by TeXtelmExtel

\begin{lemma}                         \label{simvtb}
There exists a transcendence basis $T=(x_1,\ldots,x_{\rho},y_1,\ldots,
y_{\tau}) \subset {\cal O}$ of $F|K$ with the properties:
\sn
i) \ $v_jx_i=\alpha_i\,$ for $i=1,\ldots ,\rho$, $j=1,\ldots ,\ell$,
\sn
ii) \ $v_jy_i=0$ and $y_iP_j=\xi_i$ for $i=1,\ldots ,\tau$, $j=1,\ldots
,\ell$.
\sn
Moreover, the valuations $v_1,\ldots,v_\ell$ and the places
$P_1,\ldots,P_\ell$ coincide on the rational function field
$K(T)$.
\end{lemma}
\begin{proof}
For each $i\in\{1,\ldots,\rho\}$, we employ the Approximation Theorem
for pairwise incomparable valuations ([R], Th\'eor\`eme 1, p.~135) to
find elements $x_i\in {\cal O}$ such that $v_jx_i=\alpha_i$ for
$j=1,\ldots ,\ell$. Since $\alpha_i$ does not depend on the valuation
$v_i\,$, the compatibility condition of the theorem is satisfied.

 %Inserted by TeXtelmExtel

For each $i\in\{1,\ldots,\tau\}$, we employ Theorem~11.14 of [En] to
find an element $y_i\in {\cal O}$ such that $y_iP_j=\xi_i$
for $j=1,\ldots ,\ell$.

 %Inserted by TeXtelmExtel

The last assertion is a direct consequence of Theorem~\ref{prelBour}.
\end{proof}

 %Inserted by TeXtelmExtel

\begin{remark}
The following more general assertion can be proved for pairwise
incomparable valuations with not necessarily equal value groups
(using [R], Th\'eor\`eme 3, p.~136):
{\it There are elements $z_1,\ldots,z_s\in {\cal O}$ which
simultaneously form a transcendence basis of $(F|K,P_j)$ as
described in Theorem~\ref{hrwtd}, for $j=1,\ldots ,\ell$.}
\end{remark}

 %Inserted by TeXtelmExtel

As a final ingredient for the proof of Theorem~\ref{MT3} we need a
result from the Galois theory of rings: Let $S$ be a normal, local
domain with field of fractions $E$ and maximal ideal $M$. Let $E^i(M)$
denote the absolute inertia field, i.e., the fixed field of the group
\[
I(M):=\{\sigma\in{\rm Gal} (E\sep |E):\; \forall s\in S\sep:\sigma s-s
\in M\sep\},
\]
where $S\sep$ is the integral closure of $S$ in the separable closure
$E\sep$ of $E$ and $M\sep\in\Spec S\sep$ is some prime lying over $M$.
Within $E\sep$ the field $E^i(M)$ is unique up to an $E$-isomorphism.

 %Inserted by TeXtelmExtel

Let $S^\prime$ be the integral closure of $S$ in the finite extension
$F$ of $E$ and assume that for each of the primes $N_1,\ldots ,
N_\ell\in\Spec S^\prime$ lying over $M$ the field $F$ lies in the
absolute inertia field $E^i(M)$. In particular, $F\subseteq E\sep$
holds. We consider the ring
\[
O:=\bigcap_{j=1}^\ell S^\prime_{N_j}\>;
\]
it satisfies $S^\prime\subseteq O\subseteq F$. We denote the prime
ideals $N_j S^\prime _{N_j}\cap O$ again by $N_j$. The missing
ingredient for the proof of Theorem~\ref{MT3} now reads:
\begin{lemma}
\label{multiple inertial}
If the residue field $S/M$ is infinite, then there exists a primitive
element $z\in S^\prime$ of $F|E$ with the following properties:
\begin{enumerate}
\item The minimal polynomial $p\in S[X]$ of $z$ over $E$ satisfies
$p^\prime (z)\in O^\times$.
\item The $S$-algebra $A:=S[z,p^\prime (z)^{-1}]\subseteq O$ is
standard-\'etale.
\item $A_{N_j\cap A}=O_{N_j}$ for $j=1,\ldots ,\ell$.
\end{enumerate}
\end{lemma}
\begin{proof}
We choose a prime $N_1\sep\in\Spec S\sep$ and $E$-automorphisms
$\sigma_1,\ldots ,\sigma_n$ of $E\sep$ such that their restrictions to
$F$ form a complete set of $E$-embeddings of $F$ in $E\sep$. Every
prime $N\in\Spec S^\prime$ lying over $M$ then satisfies
\begin{equation}
\label{conjugate}
N=\sigma_k^{-1}(N_1\sep\cap\sigma_k S^\prime )
\end{equation}
for some $k\in\{1,\ldots ,n\}$.

 %Inserted by TeXtelmExtel

We enumerate the $\sigma_k$ in such a way that for $j=1,\ldots ,\ell$
and  $k_0<k_1<\ldots <k_\ell$ the equations
\[
\sigma_k^{-1}(N_1\sep\cap\sigma_k S^\prime )=N_j,\; k_{j-1}\leq k<k_j,
\]
hold.

 %Inserted by TeXtelmExtel

Corresponding to equation (\ref{conjugate}) we embed the residue fields
$S^\prime /N$ over $\overline{E}:=S/M$ in $S\sep /N_1\sep$ via the
natural homomorphisms
\[
S^\prime /N\>\rightarrow\> \overline{F}_k:=\sigma_k S^\prime
/(N_1\sep\cap \sigma_k S^\prime )\>.
\]
We consider the induced ring homomorphism
\[
\phi :\; S^\prime\rightarrow\prod_{k=1}^n \overline{F}_k\>.
\]
We choose primitive elements $\vartheta_1,\ldots ,\vartheta_{k_\ell}$
for the separable extensions $\overline{F}_k|\overline{E}$, $k=1,\ldots
,k_\ell$. After multiplying with suitable elements from $\overline{E}$
we may assume the $\vartheta_k$ to be pairwise non-conjugate over
$\overline{E}$. The Chinese Remainder Theorem then provides us with an
element $z\in S^\prime$ such that $\phi z=(\vartheta_1,\ldots
,\vartheta_{k_\ell} ,0,\ldots ,0)$. This $z$ is the candidate for the
element appearing in the assertion of our lemma.

 %Inserted by TeXtelmExtel

We reduce the polynomial
\begin{equation}
\label{minpol}
p(X):=\left(\prod_{j=1}^\ell\prod_{k_{j-1}\leq k<k_j}(X-\sigma_kz)\right)
\prod_{j=k_\ell +1}^n(X-\sigma_kz)\,\in\, S\sep [X]
\end{equation}
modulo the prime $N_1\sep$. Due to the choice of $z$, we obtain
\begin{equation}
\label{redminpol}
\overline{p}(X):=\left(\prod_{j=1}^\ell\prod_{k_{j-1}\leq k<k_j}
(X-\vartheta_k)\right)X^{n-(k_\ell +1)}\,\in\, S\sep/N_1\sep [X].
\end{equation}
Since the $\vartheta_k$ are pairwise non-conjugate we conclude that
$\sigma_kz\neq z$, $k=1,\ldots ,k_\ell$, and since $\sigma_kz\in
N_1\sep$, $k=k_\ell +1,\ldots ,n$, we have $\sigma_kz\neq z$ in this
case too. Consequently, $F=E(z)$ and $p\in S[X]$ is the minimal
polynomial of $z$ over $E$. Furthermore, we conclude that $p^\prime
(\sigma_k(z))=\sigma_k (p^\prime (z))\not\in N_1\sep$ for $k=1,\ldots
,k_\ell$ and thus, $p^\prime (z)\not\in N_j$ for $j=1,\ldots ,\ell$.
This yields $p^\prime (z)\in O^\times$ and assertion 1.\ of the lemma is
verified.

 %Inserted by TeXtelmExtel

The $S$-algebra $A:=S[z,p^\prime(z)^{-1}]\subseteq O$ is $S$-isomorphic
to $(S[X]/pS[X])_{p^\prime}$, which is well-known to be
standard-\'etale.

 %Inserted by TeXtelmExtel

To verify the remaining assertion 3., we first observe that $A_{N_j\cap
A}= S[z]_{N_j\cap S[z]}$ because $p^\prime (z)\not\in N_j\cap S[z]$. On
the other hand, $A|S$ is \'etale and $S$ is normal, thus $A$ is normal
too. We conclude that $S[z]_{N_j\cap S[z]}=S^\prime_{N_j\cap S[z]}$
since this extension is integral. But then $N_j\in\Spec S^\prime$ is the
only prime lying over $N_j\cap S[z]$ and hence, $A_{N_j\cap A}=
S^\prime_{N_j\cap S[z]}=S^\prime_{N_j}$ as desired.
\end{proof}

 %Inserted by TeXtelmExtel

Finally, we turn to the

 %Inserted by TeXtelmExtel

\sn
\textbf{Proof of Theorem~\ref{MT3}}:
One starts by choosing a transcendence basis $T=\{x_1,\ldots, x_{\rho},\\
y_1,\ldots,y_{\tau}\}$ of $F|K$ according to Lemma~\ref{simvtb}.
Let $P$ denote the common restriction of $P_1,\ldots,P_\ell$ to $K(T)$.
By Theorem~\ref{hrwtd}, $(F,P_j)$ lies in the absolute inertia field of
$(K(T),P)$ for $j=1,\ldots ,\ell$. We apply Lemma~\ref{multiple
inertial} to the integral closure of ${\cal O}_{P}$ in $F$ and the ring
${\cal O}=\cap_{j=1}^\ell {\cal O}_{j}\,$. We obtain an \'etale ${\cal
O}_{P}$-algebra $A={\cal O}_{P} [z,p^\prime (z)^{-1}]\subseteq{\cal O}$.

 %Inserted by TeXtelmExtel

From now on we proceed as described in the paragraph preceding
Lemma~\ref{etale} and in its proof: We collect the structural constants of
the \'etale ${\cal O}_{P}$-algebra $A$ and the constants needed to express
the elements $\zeta\in Z$ in terms of the generators $z, p(z)^{-1}$ of $A$
into a finite set $Z^\prime\subset {\cal O}_{P}$.

 %Inserted by TeXtelmExtel

As a consequence of the assumptions made about $R$ in Theorem~\ref{MT3},
the pair $(P,Z^\prime)$ is $R$-uniformizable. Therefore there exists a
finitely generated $R$-algebra $B\subseteq K(T)$ with $\Frac B=K(T)$
such that $B_{q_B}$, $q_B:={\cal M}_{P}\cap B$, is regular and
$Z^\prime\subset B_{q_B}$. In the same way as in the proof of
Lemma~\ref{etale} we can then show that the $B_{q_B}$-algebra
$B_{q_B}[z,p^\prime (z)^{-1}]\subseteq A$ is \'etale.

 %Inserted by TeXtelmExtel

Let $C:=B[z,p^\prime (z)^{-1}]$; it follows that the localizations
$C_{{\cal M}_{j}\cap C}$, ${\cal M}_{j}$ the maximal ideal of ${\cal
O}_{j}$, are regular and contain $Z$. Hence, the morphism $\Spec
C\rightarrow\Spec B$ of affine $R$-schemes satisfies all the
requirements stated in Theorem~\ref{MT3}.
\QED

 %Inserted by TeXtelmExtel

%\bn\bn\bn
\newpage\noindent
{\bf References}
\newenvironment{reference}%
{\begin{list}{}{\setlength{\labelwidth}{5em}\setlength{\labelsep}{0em}%
\setlength{\leftmargin}{5em}\setlength{\itemsep}{-1pt}%
\setlength{\baselineskip}{3pt}}}%
{\end{list}}
\newcommand{\lit}[1]{\item[{#1}\hfill]}
\begin{reference}
\lit{[A1]} {Abhyankar, S.$\,$: {\it Uniformization in $p$-cyclic
extensions of a two dimensional regular local domain of residue field
characteristic $p$}, Festschrift zur Ged\"acht\-nis\-feier f\"ur Karl
Weierstra{\ss} 1815--1965, Wissenschaftliche Abhandlungen des Landes
Nordrhein-Westfalen {\bf 33} (1966), 243--317 (Westdeutscher Verlag,
K\"oln und Opladen)}
\lit{[A2]} {Abhyankar, S.\ S.$\,$: {\it An algorithm on polynomials in
one indeterminate with coefficients in a two dimensional regular local
domain}, Ann.\ Mat.\ Pura Appl.\ (4) {\bf 71} (1966), 25--59}
\lit{[A3]} {Abhyankar, S.$\,$: {\it Resolution of singularities of
embedded algebraic surfaces}, Academic Press, New York (1966);
2nd enlarged edition: Sprin\-ger, New York (1998)}
\lit{[A4]} {Abhyankar, S.$\,$: {\it Local uniformization on algebraic
surfaces over ground fields of characteristic $p\neq 0$}, Ann.\ Math.
\textbf{63} (1956), 491--526}
\lit{[A5]} {Abhyankar, S.$\,$: {\it Resolution of singularities of
arithmetical surfaces}, Arithmetical Algebraic Geometry, Harper and Row,
New York (1965)}
\lit{[A6]} {Abhyankar, S.$\,$: {\it Uniformization of Jungian Local
Domains}, Math.\ Annalen \textbf{159} (1965), 1--43}
\lit{[B]} {Bourbaki, N.$\,$: {\it Commutative algebra}, Hermann,
Paris (1972)}
\lit{[dJ]} {de Jong, A.~J.$\,$: {\it Smoothness, semi-stability and
alterations}, Publications math\'e\-matiques I.H.E.S.\ \textbf{83}
(1996), 51--93}
\lit{[El]} {Elliott, G.~A.$\,$: {\it On totally ordered groups, and
$K_0$}, in: Ring Theory Waterloo 1978, eds.\ D.~Handelman and
J.~Lawrence, Lecture Notes Math.\ {\bf 734}, 1--49}
\lit{[En]} {Endler, O.$\,$: {\it Valuation theory}, Berlin (1972)}
\lit{[F]} {Fu, D.$\,$: {\it Local Weak Simultaneous Resolution for High
Rational Ranks}, J.\ Algebra \textbf{194} (1996), 614--630}
\lit{[HRW]} {Heinzer, W.\ -- Rotthaus, C.\ -- Wiegand, S.$\,$: {\it
Approximating Discrete Valuation Rings by Regular Local Rings}, Proc.\
Am.\ Math.\ Soc.\ \textbf{129} No.\ 1 (2001), 37--43}
\lit{[K1]} {Kuhlmann, F.--V.$\,$: {\it Henselian function fields and tame
fields}, preprint (extended version of Ph.D.\ thesis), Heidelberg
(1990)}
\lit{[K2]} {Kuhlmann, F.--V.$\,$: {\it On local uniformization in
arbitrary characteristic}, The Fields Institute Preprint Series, Toronto
(1997)}
\lit{[K3]} {Kuhlmann, F.--V.$\,$: {\it Valuation theoretic and model
theoretic aspects of local uniformization}, in:
Resolution of Singularities - A Research Textbook in Tribute to Oscar
Zariski. Herwig Hauser, Joseph Lipman, Frans Oort, Adolfo Quiros
(eds.), Progress in Mathematics Vol.\ {\bf 181}, Birkh\"auser
Verlag Basel (2000), 381--456}
\lit{[K4]} {Kuhlmann, F.--V.$\,$: {\it On places of algebraic function
fields in arbitrary characteristic}, submitted}
\lit{[K5]} {Kuhlmann, F.--V.$\,$: {\it The generalized stability
theorem and henselian rationality of valued function fields}, in
preparation}
\lit{[K6]} {Kuhlmann, F.--V.$\,$: {\it Every place admits local
uniformization in a finite extension of the function field}, in
preparation}
\lit{[K7]} {Kuhlmann, F.--V.$\,$: {\it Valuation theory}, book in
preparation.\\
See http://math.usask.ca/$\,\tilde{ }\,$fvk/Fvkbook.htm}
\lit{[Nag]}{Nagata, M.$\,$: \textit{Finitely generated rings over a
valuation ring}, J.\ Math.\ Kyoto Univ.\ \textbf{5} (1965), 163--169}
\lit{[R]} {Ribenboim, P.$\,$: {\it Th\'eorie des valuations}, Les
Presses de l'Uni\-versit\'e de Mont\-r\'eal (1964)}
\lit{[Ray]} {Raynaud, M.$\,$: \textit{Anneaux Locaux Hens\'eliens},
Lecture Notes in Math.\ {\bf 169}, Berlin-Heidelberg-New York (1970)}
\lit{[Z]} {Zariski, O.$\,$: {\it Local uniformization on
algebraic varieties}, Ann.\ Math.\ {\bf 41} (1940), 852--896}
\lit{[Z2]} {Zariski, O.$\,$: {\it The reduction of singularities of an
algebraic surface}, Ann.\ Math.\ {\bf 40} (1939), 639--689}
\lit{[Z3]} {Zariski, O.$\,$: {\it A simplified proof for resolution of
singularities of an algebraic surface}, Ann.\ Math.\ {\bf 43} (1942),
583--593}
\end{reference}
\adresse
\end{document}